\documentclass[11pt]{article}
\input amssymb.sty
\input amssym.def
\input amssym
\input epsf
\newtheorem{theorem}{Theorem}[section]

\newtheorem{theorem-definition}[theorem]{Theorem-Definition}
\newtheorem{theorem-construction}[theorem]{Theorem-Construction}
\newtheorem{lemma-definition}[theorem]{Lemma--Definition}
\newtheorem{proposition-definition}[theorem]{Proposition--Definition}
\newtheorem{lemma}[theorem]{Lemma}
\newtheorem{proposition}[theorem]{Proposition}
\newtheorem{corollary}[theorem]{Corollary}

\newtheorem{definition}[theorem]{Definition}

\begin{document}
\newcommand{\Z}{{\mathbb Z}}
\newcommand{\R}{{\mathbb R}}
\newcommand{\Q}{{\mathbb Q}}
\newcommand{\C}{{\mathbb C}}
\newcommand{\lms}{\longmapsto}
\newcommand{\lra}{\longrightarrow}
\newcommand{\hra}{\hookrightarrow}
\newcommand{\ra}{\rightarrow}
\newcommand{\sgn}{\rm sgn}
\begin{titlepage}\title{Pentagon relation for the quantum dilogarithm and quantized
${\cal M}^{\rm cyc}_{0,5}$}
\author{A.B. Goncharov}
\end{titlepage}
\date{\it To the memory of Sasha Reznikov}
%\stepcounter{page}
\maketitle
%\addtocounter{page}{+1}
%\doublespace
\tableofcontents

 %vskip 6mm \noindent

\section{Introduction}

Let $\hbar>0$. The quantum dilogarithm
function is given by the following integral:
$$
\Phi^\hbar(z) := {\rm exp}\Bigl(-\frac{1}{4}\int_{\Omega}\frac{e^{-ipz}}{ {\rm sh} (\pi p)
{\rm sh} (\pi \hbar p) } \frac{dp}{p} \Bigr),
\qquad {\rm sh}(p)= \frac{e^p-e^{-p}}{2}.
$$
Here $\Omega$ is a path from $-\infty$ to $+\infty$ making a little
half circle going over the zero.
So the integral is convergent.
It goes back to  Barnes \cite{Ba}, and appeared in many papers during the last 30 years:
\cite{Bax}, \cite{Sh}, \cite{Fad1}, ... .
The function $\Phi^\hbar(z)$  enjoys
 the following properties (cf. \cite{FG3}, Section 4):

\begin{itemize}

\item   The function $\Phi^\hbar(z)$ is meromorphic. Its zeros are simple zeros
in the upper half plane at the points
\begin{equation}\label{tu1}
\{\pi i ((2m-1)+
(2n-1)\hbar)|m,n \in {\mathbb N}\},\qquad {\mathbb N}:= \{1, 2, \ldots\}.
\end{equation}
  Its poles are simple poles, located
in the lower half plane, at the
points
\begin{equation} \label{tu2}
\{-\pi i ((2m-1)+(2n-1)\hbar)|m,n \in {\mathbb N}\}.
\end{equation}

\item  The function $\Phi^\hbar(z)$  is characterized by the following difference relations. 
Let $q:= e^{\pi i \hbar}$ and  $q^{\vee}:= e^{\pi i /\hbar}$. Then 
\begin{equation}\label{tu}
\Phi^\hbar(z+ 2 \pi i h) = \Phi^\hbar(z) (1+qe^z), \quad
\Phi^\hbar(z+ 2 \pi i ) = \Phi^\hbar(z) (1+ q^{\vee}e^{z/\hbar}), 
\end{equation}

\item  One has $|\Phi^\hbar(z)|=1$ when $z$ is on the real line.

\item  It is related in several ways to the dilogarithm, e.g.
its asymptotic expansion when $\hbar \to 0$ is
$$
\Phi^\hbar(z) \sim {\rm exp}\Bigl(\frac{{\rm L}_2(e^{z})}{2\pi i \hbar}\Bigr),
\quad \mbox{where } {\rm L}_2(x):= \int_0^x\log(1+t)\frac{dt}{t}
$$
is a version of the Euler's dilogarithm function.
\end{itemize}

When $\hbar$ is
a complex number with ${\rm Im}~\hbar > 0$,  there is an infinite product expansion
$$
\Phi^\hbar(z) = \frac{{\Psi}^q(e^{z})}
{{\Psi}^{1/q^{\vee}}(e^{z/\hbar})}, \qquad \mbox{where}  ~~~
{\bf \Psi}^q(x):= \prod_{a=1}^{\infty}(1+q^{2a-1}x)^{-1}.
$$

The function $\Phi^\hbar(z)$ provides an operator $K: L^2(\R) \to L^2(\R)$,
defined as a rescaled
Fourier transform followed by
the operator of multiplication by the quantum dilogarithm $\Phi^\hbar(x)$.
$$
Kf(z):= \int_{-\infty}^{\infty}f(x)\Phi^\hbar(x) {\rm exp}(\frac{-xz}{2\pi i \hbar})dx.
$$
Since $|\Phi^\hbar(x)|=1$ on the real line, $2\pi \sqrt {\hbar}K$ is unitary.

\begin{theorem} \label{qp4*}
$(2\pi  \sqrt {\hbar} K)^5 = \lambda \cdot {\rm Id}$, where
$|\lambda|= 1$.
\end{theorem}
In the quasiclassical limit it gives Abel's five term
relation for the dilogarithm.

\vskip 3mm
The  pentagon relation for the simpler version
${\bf \Psi}^q(x)$
of the quantum dilogarithm was discovered in \cite{FK}. A similar pentagon relation
for the function $\Phi^\hbar(z)$, which is
equivalent to Theorem \ref{qp4*},  was suggested in \cite{Fad1} and proved,
using different methods, in \cite{Wo} and
\cite{FKV}.
Theorem \ref{qp4*} was formulated in \cite{CF}. However
the argument presented there as a proof has a significant problem,
which put on hold the program of quantization of
Teichm\"uller spaces.

\vskip 3mm
In this paper we show that the operator $K$ is a part of a much more rigid structure, called
the {\it quantized moduli space ${\cal M}^{\rm cyc}_{0,5}$} -- this easily implies
 Theorem \ref{qp4*}.

\begin{figure}[ht]
\centerline{\epsfbox{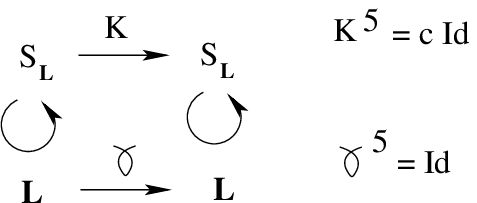}}
\caption{Quantized moduli space ${\cal M}^{\rm cyc}_{0,5}$.}
\label{qp2a}
\end{figure}

Namely, consider the algebra
generated by operators of multiplication by
$e^x$ and $e^{x/\hbar}$ and shifts by $2\pi i$ and $2\pi i\hbar$, acting as unbounded operators in
$L^2(\R)$.
We use  a remarkable subalgebra ${\bf L}$ of
this $\ast$-algebra, and
introduce a {\it Schwartz space $S_{\bf L} \subset L^2(\R)$}, defined as the common domain of the
operators from ${\bf L}$. It comes with a natural topology.
Our main result, Theorem \ref{qp21}, tells that the operator $K$ preserves the space $S_{\bf L}$, and the
conjugation by $K$ intertwines an order
$5$ automorphism $\gamma$ of the algebra ${\bf L}$, see Fig. \ref{qp2a}.
This characterises the operator $K$ up to a constant.
The proof uses  analytic properties of the space ${S}_{\bf L}$ developed in
Theorem \ref{qp12}. Theorem \ref{qp21} easily implies Theorem \ref{qp4*}.

We define a space of {\it distributions}
$S^*_{\bf L}$ as the topological dual to ${S}_{\bf L}$.
So there is a Gelfand triple $S_{\bf L}\subset L_2(\R) \subset S^*_{\bf L}$. The  operator $K$ acts by
its automorphisms.
It would be interesting to calculate it on some distributions explicitly.

\vskip 3mm
The story is similar in spirit to the Fourier
 transform theory developed using the algebra
of polynomial differential operators:
$$
 \mbox{The Fourier transform} ~<-> ~\mbox{The operator $K$}.
$$
$$
 \mbox{The algebra ${\bf D}$ of polynomial differential operators}~ <->~
\mbox{The algebra ${\bf L}$ of difference operators}.
$$
$$
\mbox{The automorphism $\varphi$ of  ${\bf D}$ given by $ix \to d/dx$, $d/dx \to -ix$}
~<-> ~\mbox{The automorphism $\gamma$ of ${\bf L}$}.
$$
$$
\mbox{The classical Schwartz space} ~<-> ~\mbox{The Schwartz space ${S}_{\bf L}$}.
$$
\vskip 3mm
Let $  {\cal M}^{\rm cyc}_{0,5}\subset \overline {\cal M}_{0,5}$ be the moduli space
of configurations of $5$ cyclically ordered points on ${\Bbb P}^1$, where we do not allow
the neighbors to collide.
It carries an atlas consisting of $5$ coordinate systems,
providing $  {\cal M}^{\rm cyc}_{0,5}$ with a structure of the 
cluster ${\cal X}$-variety of type $A_2$. The algebra ${\bf L}$ is isomorphic to the
algebra of regular functions on the modular double of  the non-commutative $q$-deformation
of the cluster ${\cal X}$-variety.
The automorphism $\gamma$ corresponds to a cyclic shift acting on
configurations of points.
\vskip 3mm
{\it The triple $({\bf L}, {S}_{\bf L}, \gamma)$, see Fig. \ref{qp2a}, is called  the
  quantized moduli space
$ {\cal M}^{\rm cyc}_{0,5}$}.

\vskip 3mm
The results of this paper admit a
generalization to a cluster set-up,
where the role of the automorphism $\gamma$ plays the cluster mapping class group.
In particular this gives a definition of quantized higher Teichm\"uller spaces, and allows to
state  precisely the modular functor property of the latter.
\footnote{In previous versions %(\cite{K}, \cite{CF}, \cite{FG2}-\cite{FGII})
of
quantization of Teichm\"uller spaces/cluster ${\cal X}$-varieties
 the pair (${\bf L}$,
$S_{\bf L}$) was missing, making the resulting notion rather flabby.}

\vskip 3mm
{\it The structure of the paper}. In Section 2.1 we recall
the cluster ${\cal X}$-variety of type $A_2$ \cite{FG2}.
In Section 3 we identify it with 
$  {\cal M}^{\rm cyc}_{0,5}$.
This clarifies formulas in Section 2.1-2.2.
 In Section 2.2 we recall 
a collection of regular functions on our cluster ${\cal X}$-variety.
Theorem \ref{pq2} tells that they form a basis in the space of regular functions,
and in particular closed under multiplication.
We introduce a $q$-deformed version of this basis/algebra. Its tensor product
with a similar algebra for $q^{\vee}$ is the algebra ${\bf L}$.
In Sections 2.3-2.5 we prove our main results.

\vskip 3mm
{\bf Acknowledgments}.
I am very grateful to Joseph Bernstein  for several illuminating discussions,
 and to Andrey Levin for reading carefully the first draft of the text and spotting some errors.

I was
supported by the  NSF grants  DMS-0400449 and DMS-0653721.
The final version of the text was written during my stay in IHES.
I would like to thank IHES for hospitality and support.

\section{Quantized moduli space $  {\cal M}^{\rm cyc}_{0,5}$}

\subsection
{Cluster varieties of type $A_2$}
The cluster ${\cal X}$-variety is glued from five copies of
$\C^* \times \C^*$, so that $i$-th copy is glued to $(i+1)$-st
(indexes are modulo $5$)
by the map acting on the coordinate functions as follows:\footnote{We use a definition which differs
slightly from the standard one, but delivers the same object.}
\begin{equation} \label{1}
\gamma^*_X: X \lms Y^{-1}, \quad Y \lms (1+Y)X.
 \end{equation}
Similarly, the cluster ${\cal A}$-variety is glued from five copies of
$\C^* \times \C^*$, so that $i$-th copy is glued to $(i+1)$-th by the map
acting on the coordinate functions as follows:
\begin{equation} \label{2}
\gamma^*_A: (A, B)\mapsto
((1+A)B^{-1}, A).
\end{equation}
(Accidently, these two cluster varieties are canonically isomorphic).

The fifth degree of each of these maps is the identity. Thus
the map identifying the $i$-th copy of $\C^* \times \C^*$ with the $(i+1)$-st one
in the standard way is an automorphism of order $5$
acting on the  ${\cal X}$- and ${\cal A}$-varieties.
We denote it by $\gamma$.

Recall the tropical semifield $\Z^t$. It is the set $\Z$ with the operations of addition
$a \oplus b:= {\rm max}\{a,b\}$, and multiplication $a\otimes b := a+b$.
The set ${\cal A}(\Z^t)$ of $\Z^t$-points of the  ${\cal A}$-variety is defined by gluing
the five copies of $\Z^2$ via the tropicalizations of the map (\ref{2}).
The map  $\gamma$ acts
on the tropical ${\cal A}$-space by
$$
\gamma_a: (a,b)\mapsto
(\max(a,0)-b, a),\quad \gamma^5={\rm Id}.
$$

\begin{figure}[ht]
\centerline{\epsfbox{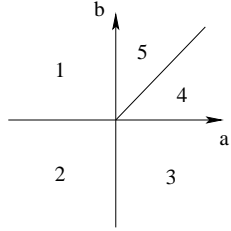}}
\caption{The five domains in the tropical ${\cal A}$-space.}
\label{qp1}
\end{figure}

There are five cones in the
tropical ${\cal A}$-space, shown on Fig. \ref{qp1}.
The map $\gamma$ shifts them cyclically counterclockwise.
It is a piecewise linear map, whose restriction to each cone  is linear.

\subsection{The $\ast$-algebra {\bf L}}

{\bf The canonical
basis  for the cluster ${\cal X}$-variety  of type $A_2$.}
A rational function $F(X,Y)$ is a {\it
universally positive Laurent polynomial} on ${\cal X}$, if
$(\gamma_X^*)^iF(X,Y)$ is a Laurent polynomial with positive integral coefficients
for every $i$.
Equivalently, it belongs to  the intersection of the ring of regular functions
on the scheme ${\cal X}$ over $\Z$ with the semifield of
rational functions with positive integral coefficients.
There is a canonical $\gamma$-equivariant map, defined in Section 4 of \cite{FG2}:
\footnote{Observe that $\gamma^*_X$ tells how the automorphism $\gamma$ acts on functions, while $\gamma_a$ tells the action on the tropical points.}
$$
{\Bbb I}_{\cal A}: {\cal A}(\Z^t) \lra
\mbox{The space of universally positive Laurent polynomials on ${\cal X}$},
$$
\begin{equation} \label{pop}
\gamma_X^{*}{\Bbb I}_{\cal A}(\gamma_a(a,b)) = {\Bbb I}_{\cal A}(a,b),
\end{equation}
given by:
\begin{equation} \label{rt}
{\Bbb I}_{\cal A}(a,b)=\left\{\begin{array}{lll} X^{a}Y^b& \mbox{ for } & a\leq 0
\mbox{ and
} b\geq 0\\
\left(\frac{1+X}{XY}\right)^{-b}X^{a}& \mbox{ for } & a\leq 0
\mbox{
and } b\leq 0\\
\left(\frac{1+X+XY}{Y}\right)^a\left(\frac{1+X}{XY}\right)^{-b}&
\mbox{ for } & a\geq 0 \mbox{ and } b\leq 0\\
((1+Y)X)^{b}\left(\frac{1+X+XY}{Y}\right)^{a-b}&
\mbox{ for } & a\geq b\geq 0\\
Y^{b-a}((1+Y)X)^a &\mbox{ for } & b\geq a\geq 0.\\
\end{array}
\right.
\end{equation}
Or equivalently, showing that the leading monomial is always $X^aY^b$: 
$$
{\Bbb I}_{\cal A}(a,b)=\left\{\begin{array}{lll} X^{a}Y^b& \mbox{ for } & a\leq 0
\mbox{ and } b\geq 0\\
X^a Y^b(1+X^{-1})^{-b} & \mbox{ for } & a\leq 0 \mbox{
and } b\leq 0\\
X^a Y^b(1+X^{-1})^{-b}(1+Y^{-1}+X^{-1}Y^{-1})^a&
\mbox{ for } & a\geq 0 \mbox{ and } b\leq 0\\
X^a Y^b (1+Y^{-1})^b(1+Y^{-1}+X^{-1}Y^{-1})^{a-b}&
\mbox{ for } & a\geq b\geq 0\\
X^aY^b(1+Y^{-1})^a&\mbox{ for } & b\geq a\geq 0.\\
\end{array}
\right.$$ 
One can easily verifies that the formulae agree on the
overlapping domains of values of $a$ and $b$.
The $i$-th row of (\ref{rt})
describes the restriction of the canonical map to the $i$-th cone.

\vskip 3mm

{\bf The quantum ${\cal X}$-variety
 and the quantum canonical basis}.
Let ${\rm T}_q$ be the algebra generated over $\Z[q,q^{-1}]$ by
$X^{\pm 1}, Y^{\pm 1}$, subject to the relation
$q^{-1}XY -qYX = 0$. It is called the two dimensional quantum torus algebra.
It has an involutive antiautomorphism $\ast$ such that
$$
\ast q = q^{-1}, \quad \ast X = X, \quad \ast Y = Y.
$$
Consider the following $q$-deformation of the $\ast$-equivariant map $\gamma$:
\begin{equation} \label{qp7}
\gamma^*_q: X \lms Y^{-1}, \quad Y \lms (1+qY)X.
\end{equation}
One checks that it  is an order $5$ automorphism of the fraction field of
${\rm T}_q$.  The quantum ${\cal X}$-space ${\cal X}_q$ is nothing else but
a pair $({\rm T}_q, \gamma^*_q)$.
\footnote{Alternatively,
using a geometric language, the quantum ${\cal X}$-space ${\cal X}_q$
is glued from
five copies of the spectrum ${\rm Spec}({\rm T}_q)$ of the quantum torus ${\rm T}_q$,
so that $i$-th copy is glued to $(i+1)$-st
along the map (\ref{qp7})}.

\vskip 3mm
An element $F(X,Y)$ of the fraction field of ${\rm T}_q$ is a {\it
universally positive Laurent polynomial} on ${\cal X}_q$ if
$(\gamma_q^*)^iF(X,Y)$ is a Laurent polynomial in $X,Y,q$
with positive integral coefficients
for every $i$.

\begin{proposition} \label{qp6}
There is a canonical $\gamma$-equivariant map
$$
{\Bbb I}^q_{\cal A}: {\cal A}(\Z^t) \lra
\mbox{The space of universally positive Laurent polynomials on ${\cal X}_q$}.
$$
\end{proposition}

\vskip 3mm
{\bf Construction}.
It is obtained by multiplying each monomial in (\ref{rt})
by a (uniquely defined) power of $q$, making it $\ast$-invariant. For example,
the quantum canonical map on the first cone is given by
$$
{\Bbb I}^q_{\cal A}(a,b) = q^{-ab}X^aY^b, \quad a\leq 0, b\geq 0.
$$
Then we can use (\ref{pop}), which is valid in the $q$-deformed version as well. 
%\begin{equation} \label{rt1}
%{\Bbb I}_{\cal A}(a,b)=\left\{\begin{array}{lll} q^{-ab}X^{a}Y^b& \mbox{ for } & a\leq 0
%\mbox{ and
%} b\geq 0\\
%\left(\frac{1+X}{XY}\right)^{-b}X^{a}& \mbox{ for } & a\leq 0
%\mbox{
%and } b\leq 0\\
%\left(Y^{-1}+qXY^{-1}+X\right)^a\left(\frac{1+X}{XY}\right)^{-b}&
%\mbox{ for } & a\geq 0 \mbox{ and } b\leq 0\\
%((1+qY)X)^{b}\left(Y^{-1}+qXY^{-1}+X\right)^{a-b}&
%\mbox{ for } & a\geq b\geq 0\\
%Y^{b-a}((1+qY)X)^a &\mbox{ for } & b\geq a\geq 0.\\
%\end{array}
%\right.
%\end{equation}
\footnote{We do not use the
fact that  Laurent $q$-polynomials in the basis
have positive integral coefficients.}

\vskip 3mm
Denote by ${\Bbb L}_q$
the image of the map ${\Bbb I}^q_{\cal A}$. It is
closed under multiplication.
Set ${\bf L}:= {\Bbb L}_q \otimes {\Bbb L}_{q^{\vee}}$.

\vskip 3mm
\subsection{The Schwartz space  $S_{\bf L}$}
Let $W \subset L^2(\R)$ be the space of finite
$\C$-linear combinations of the functions
\begin{equation} \label{7.18.07.1}
e^{-ax^2/2+bx}P(x), \quad \mbox{
where $P(x)$ is a polynomial in $x$, and $a\in \R_{>0}, b \in \C$}.
\end{equation}
Set
$$
\widehat X (f)(x):= f(x+2\pi i \hbar), \quad \hbar \in \R_{>0},
\qquad \widehat Y (f)(x):= e^xf(x).
$$
$$
\widehat X^{\vee} (f)(x):= f(x+2\pi i), \quad \hbar \in \R_{>0},
\qquad \widehat Y^{\vee} (f)(x):= e^{x/\hbar}f(x).
$$
They are symmetric unbounded operators.
They preserve $W$
and satisfy, on $W$,   relations
$$
\widehat X \widehat Y  = q^{2}\widehat Y\widehat X, \quad q:= e^{\pi i \hbar}.
$$
$$
\widehat X^{\vee} \widehat Y^{\vee}  =
(q^{\vee})^{2}\widehat Y^{\vee}\widehat X^{\vee}, \quad q^{\vee}:= e^{\pi i /\hbar}.
$$
The second pair of operators commute with the first one.
Therefore these  operators provide an $\ast$-representation of
the algebra ${\rm T}_q \otimes {\rm T}_{q^{\vee}}$ in $W$.

\vskip 3mm
{\bf Remark}. Consider a smaller subspace $W_0 \subset W$,
with
$a=1, b \in 2\pi i \hbar \Z + 2 \pi i \Z + \Z + 1/\hbar\Z$  and
${\rm deg}(P) =0$. Then acting on $e^{-x^2/2}$ we get an
isomorphism of linear spaces
${\rm T}_q \otimes {\rm T}_{q^{\vee}} \stackrel{\sim}{=} W_0$.

\vskip 3mm
In particular an  element $A \in {\bf L}$
 acts by an
 unbounded
operator $\widehat A$ in $W$.

\begin{definition} \label{7.10.07.1} The Schwartz space $S_{{\bf L}}$
for the $\ast$-algebra  ${\bf L}$
is a subspace of $L^2(\R)$ consisting
of vectors $f$ such that the functional $w \to (f, \widehat Aw)$ on $W$
is continuous for the $L_2$-norm.
\end{definition}
Denote by $(*,*)$ the scalar product in $L_2$.
The Schwartz space for the $\ast$-algebra ${\bf L}$
is
 the common domain of definition of
operators from ${\bf L}$ in $L^2(\R)$. Indeed,
since $W$ is dense in  $L^2(\R)$,
the Riesz theorem implies that for any $f \in S_{\bf L}$ there exists a
unique $g \in L^2(\R)$ such that $(g,w) = (f, \widehat Aw)$. We set $\widehat {\ast A}f := g$.
Equivalently,
let $W^*$ be the algebraic linear dual to $W$.
So $L^2(\R) \subset W^*$. Then
$$
S_{\bf L}= \{v \in W^*| \widehat A^*v \in L^2(\R) \quad \mbox{for any $A \in {\bf L}$}\} \cap
L^2(\R).
$$

The Schwartz space $S_{{\bf L}}$ has a natural topology given by  seminorms
$$
\rho_B(f):= ||Bf||_{L_2}, \qquad \mbox{$B$ runs through a basis in ${\bf L}$}.
$$

The key properties of the Schwartz space $S_{{\bf L}}$ which we use below are the following.
\begin{theorem} \label{qp12}
The space $W$ is dense in the Schwartz space $S_{{\bf L}}$.
\end{theorem}

One can interpret 
Theorem \ref{qp12} by saying that
the {\it $\ast$-algebra ${\bf L}$ is  essentially self-adjoint in $L^2(\R)$}.

\vskip 3mm

{\bf Proof}. %(i).

\begin{lemma} \label{qp141}
For any $w \in W, s \in S_{{\bf L}}$, the convolution
$s \ast w $ lies in $S_{{\bf L}}$.
\end{lemma}

{\bf Proof}. Set
 $T_\lambda f(x) := f(x-\lambda)$. Write
$$
s \ast w (x)= \int_{-\infty}^{\infty}w(t) (T_ts)(x)dt.
$$
For any seminorm $\rho_B$ on $S_{{\bf L}}$ the operator
$T_\lambda: (S_{{\bf L}}, \rho_B) \lra (S_{{\bf L}}, \rho_B)$ is a bounded operator
with the norm bounded by $e^{|\lambda|}$. Thus the operator
$\int_{-\infty}^{\infty}w(t) T_tdt$ is a bounded operator on $(S_{{\bf L}}, \rho_B)$.
This implies the Lemma.

\vskip 3mm
 Let $w_{\varepsilon}:= (2\pi)^{-\frac{1}{2}}\varepsilon^{-1} e^{-\frac{1}{2}(x/\varepsilon)^2}\in W$ be a sequence
converging as $\varepsilon\to 0$ to  the $\delta$-function at $0$.
Clearly one has in the topology of
$S_{{\bf L}}$
\begin{equation} \label{qp13}
\lim_{\varepsilon\to 0}w_{\varepsilon} \ast s  = s(x).
\end{equation}

\begin{lemma} \label{qp14}
For any $w \in W, s \in S_{{\bf L}}$,
the Riemann sums for the integral
\begin{equation} \label{qp131}
s \ast w (x) =
\int_{-\infty}^{\infty}s(t) \omega(x-t)dt = \int_{-\infty}^{\infty}s(t) T_tw(x)dt.
\end{equation}
 converges in the topology of
$S_{{\bf L}}$ to the convolution $s \ast w$.
\end{lemma}

{\bf Proof}. Let us show first that (\ref{qp131})  is convergent in
$L_2(\R)$. The key fact is that a shift of
$w\in W$  quickly becames essentially orthogonal to $w$. More precisely,
in the important for us case when
 $w = {\rm exp}(-ax^2/2+bx)$, $a>0$,  (this includes any $w \in W_0$)  we have
\begin{equation} \label{qp132}
( w(x), T_\lambda w(x)) < C_we^{-a\lambda^2/2+ (b-\overline b) \lambda}.
\end{equation}
Therefore in this case
$$
\Bigl(\int_{-\infty}^{\infty}s(t) T_tw(x)dt, \int_{-\infty}^{\infty}s(t) T_tw(x)dt\Bigr)
$$
$$
\stackrel{(\ref{qp132})}{\leq} C_w\int_{-\infty}^{\infty}\int_{-\infty}^{\infty}
e^{-a(t_1-t_2)^2/2 + (b-\overline b)(t_1-t_2)} |s(t_1) s(t_2)| dt_1 dt_2 =
$$
$$
C_w\int_{-\infty}^{\infty}e^{-a \lambda^2/2+(b-\overline b)\lambda}
\int_{-\infty}^{\infty}|s(t) s(t+\lambda)|  dtd\lambda \leq
 C_w||s||_{L_2}^2\int_{-\infty}^{\infty}e^{-a\lambda^2/2+(b-\overline b)\lambda}  d\lambda.
$$
We leave the case of an arbitrary $w$ to the reader: it is not used in the proof of the theorem.

The convergence with respect to the seminorm $||Bf||$ is proved by the same argument using
the fact that $W_0 \subset W$ is stable under the algebra ${\bf L}$.
The Theorem follows from Lemma \ref{qp141},  (\ref{qp13})  and Lemma \ref{qp14}.

\vskip 3mm
{\bf Remark}.  The same arguments show that the space $W$ is dense
in the space $S_{\bf L'}$ defined for any subalgebra ${\bf L'}$ of ${\bf L}$.

\vskip 3mm
\subsection{The main result}
\begin{theorem} \label{qp21}
The operator $K$ preserves the Schwartz space $S_{\bf L}$.
It
intertwines the automorphism $\gamma$ of  the algebra ${\bf L}$, i.e.
for any $A \in {\bf L}$ and $s \in S_{\bf L}$ one has
\begin{equation} \label{tu31}
K^{-1}\widehat A K s = \widehat {\gamma(A)} s.
\end{equation}
\end{theorem}

{\bf Proof}. We need the following key result.

\begin{proposition} \label{qp2}
For any $A \in {\bf L}$, $w \in W$ one has
$
K\widehat {\gamma(A)}w = \widehat AKw.
$ Therefore $\widehat AKw \in L^2(\R)$.
\end{proposition}

{\bf Proof}.
Let ${\Bbb L}_q'$ be the space of Laurent $q$-polynomials $F$ in $X, Y$
such that $\gamma(F)$ is again  a Laurent $q$-polynomial.
The following elements belong to ${\Bbb L}_q'$:
\begin{equation} \label{qp1001}
X^aY^m, \qquad
X^aY^{-n}(1+qX^{-1}) (1+q^3X^{-1}) \ldots(1+q^{2n-1}X^{-1}) , \qquad a \in \Z, \quad m,n \geq 0.
\end{equation}
Indeed, $
\gamma(Y^{-n}) = ((1+qY)X)^{-n} = X^{-n}\prod_{a=1}^n(1+q^{(2a-1)}Y)^{-1}.
$

\begin{lemma} \label{qp1000}    (i) The monomials (\ref{qp1001}) span the space ${\Bbb L}_q'$.

(ii) For every $A \in {\Bbb L}_q'\otimes {\Bbb L}'_{q^{\vee}}$, $w \in W$ one has $
K\widehat {\gamma(A)}w = \widehat AKw.
$
\end{lemma}

{\it Lemma \ref{qp1000} implies Proposition \ref{qp2}}. Thanks to
the very definition ${\Bbb L}_q \subset {\Bbb L}'_q$. So the part (ii) of the Lemma 
includes 
the commutation relation from Proposition \ref{qp2}.

\vskip 3mm
{\bf Proof of Lemma \ref{qp1000}}. (i) is obvious.

(ii) Let us prove first the following three basic identities:
\begin{equation} \label{tu3}
K (1+q\widehat Y)\widehat X w = \widehat Y Kw; \qquad
K \widehat Y^{-1}  w = \widehat XKw; \quad K \widehat X^{-1} w =
\widehat Y^{-1}(1+q\widehat X^{-1})   Kw.
\end{equation}
The general case follows from this. To see this, observe that if $A_1, A_2 \in {\bf L}$ and 
$K\widehat A_i w = \widehat {\gamma^{-1}(A_i)}Kw$ for $i=1,2$, then, since  $\widehat A_2 w \in W$, one has 
$$
K \widehat A_1 \widehat A_2 w = \widehat {\gamma^{-1}(A_1)} K \widehat A_2 w = 
\widehat {\gamma^{-1}(A_1A_2)}Kw. 
$$

{\it The first identity}. Denote by $C_s$ the line $x+is$ parallel to the $x$-axis. One has
$$
K (1+q\widehat Y)\widehat X w = \int_{C_0}(1+qe^{x})w(x+2\pi i \hbar )
\Phi^\hbar(x) e^{-xz/2\pi i \hbar}dx
$$
$$
\stackrel{(\ref{tu})}{=} \int_{C_0}w(x+2\pi i \hbar )
\Phi^\hbar(x+2\pi i \hbar ) {e}^{-xz/2\pi i h}dx
= \int_{C_{2\pi i\hbar}}w(x)
\Phi^\hbar(x) {e}^{-(x-2\pi i\hbar)z/2\pi i h}dx
$$
$$
\stackrel{(\ast)}{=} e^{z}Kw=
\widehat YKw.
$$
To obtain the equality $(\ast)$  we have to move the contour
$C_{2\pi i \hbar}$ down to $C_{0}$.
We can justify this because: (i) the function $\Phi^\hbar(z)$ is analytic
in the upper half plane, see (\ref{tu2}), and (ii) the function $\Phi^\hbar(z)$ growth
on any horizontal line at most exponentially, while $w(x)$ decays there much faster, like $e^{-x^2}$.

\vskip 3mm
{\bf Remark}.  We used here $\hbar >0$. We would not be able to
move a similar contour with negative imaginary part, since  it will hit the poles of
$\Phi^\hbar(z)$.
\vskip 3mm

{\it The second identity}.
$$
K \widehat Y^{-1}  w(z) = K e^{-x} w = \int_{-\infty}^{\infty}w(x)
\Phi^\hbar(x) {e}^{-x(z+2\pi i \hbar)/2\pi i \hbar}dx =
\widehat XKw(z).
$$

{\it The third identity}.
We have
$$
\int_{C_0} w(x-2\pi i \hbar) \Phi^\hbar(x) e^{-xz/2\pi i \hbar}dx 
$$
$$
=\int_{C_{-2\pi i \hbar}} w(x) \Phi^\hbar(x+2\pi i \hbar) e^{-(x+2\pi i \hbar) z/2\pi i \hbar}dx 
$$
$$
=e^{-nz} \int_{C_{-2\pi i \hbar}} w(x) \Phi^\hbar(x+2\pi i \hbar) e^{-x z/2\pi i \hbar}dx.
$$
We can move the contour $C_{-2\pi i \hbar}$ up towards $C_0$ since the
function $\Phi^\hbar(x+2\pi i \hbar)$ is holomorphic above the line $C_{-2\pi i \hbar}$,
and grows in horizontal directions
in the area between the two contours at most exponentially, while $w(x)$ decays like $e^{-x^2}$.
 So we get
$$
e^{-z} \int_{C_{0}} w(x) \Phi^\hbar(x+2\pi i \hbar) e^{-x z/2\pi i \hbar}dx 
$$
$$
=e^{-z} \int_{C_{0}} w(x) \Phi^\hbar(x) (1+qe^x)  e^{-x z/2\pi i \hbar}dx
$$
$$
= Y^{-1} (1+q\widehat X^{-1})  \int_{C_{0}} w(x) \Phi^\hbar(x) e^{-x z/2\pi i \hbar}dx.
$$
Lemma \ref{qp1000} is proved.

\vskip 3mm

To show that  $Ks \in S_{\bf L}$ for an $s \in  S_{\bf L}$ we need to check that
for any $B \in {\bf L}$ the functional $w \to (Ks, \widehat B^*w)$ is continuous.
Since $W$ is dense in $S_{\bf L}$ by Theorem \ref{qp12}, there is a sequence $v_i \in W$ converging to
$s$ in  $S_{\bf L}$. This means that
\begin{equation} \label{tu32}
\lim_{i\to \infty}(\widehat Bv_i, w) = (\widehat Bs,w)\quad \mbox{for any $B \in {\bf L}$, $w \in W$}.
\end{equation}
One has
$$
(Ks, \widehat B^*w) =
(s, K^{-1}\widehat B^*w)
\stackrel{\mbox{Prop. \ref{qp2}}}{=}
(s, \widehat {\gamma(B^*)}K^{-1}w) \stackrel{(\ref{tu32})}{=}
\lim_{i\to \infty} (v_i, \widehat {\gamma(B^*)}K^{-1}w) 
$$
$$
\stackrel{\rm def}{=} 
\lim_{i\to \infty} (\widehat {\gamma(B)}v_i, K^{-1}w) \stackrel{(\ref{tu32})}{=}
(\widehat {\gamma(B)}s, K^{-1}w) =  (K\widehat {\gamma(B)}s, w).
$$
Since the functional on the right is continuous,  $Ks \in S_{\bf L}$, and we have (\ref{tu31}).
The theorem is proved.
\vskip 3mm
Since $\gamma^5 = {\rm Id}$, Theorem \ref{qp21} immediately implies
\begin{corollary} \label{qp5}
For any $A \in {\bf L}$ one has
$K^{-5}\widehat AK^5 = \widehat A $ on $S$.
\end{corollary}

\subsection{Proof of Theorem \ref{qp4*}}
Let
$
E = \{f\in L^2(\R)| e^{nx}f(x) \in L^2(\R) \quad \mbox{for any $n>0$}\}.
$
\begin{lemma} \label{p1}
$
K^5(E) \subset E.
$
\end{lemma}

{\bf Proof}. Indeed, since $\widehat Y = e^x$ and $Y^n \in {\bf L}$ for any $n>0$,
one has $K^5e^{nx}f = e^{nx}K^5f$ for any $n>0$, $f \in S$
by Corollary \ref{qp5}. So using the Remark in the end of Section 2.3, we see that $W$, and hence
$S$ is dense in $E$, we get the claim.

\vskip 3mm
\begin{lemma}\label{per3}
$K^5$ is the operator of multiplication by a function
$F(x)$.
\end{lemma}

{\bf Proof}. We claim that the value $(K^5f)(a)$ depends only on the value
$f(a)$. Indeed, for any
$f_0(x) \subset E, f_0(a) = f(a)$ we have
 $f = (e^x-e^a)\phi(x) + f_0(x)$, where
$
\phi(x)= (f-f_0)/(e^x-e^a)\in E.
$
Thus
$
K^5f = (e^x-e^a)K^5\phi(x) + K^5f_0(x).
$
So $K^5f(a) = K^5f_0(a)$. Now define  $F(a)$ from
$K^5f_0(a) = F(a)f_0(a)$. The lemma is proved.

\vskip 3mm
\begin{proposition} \label{per1}
The function $F(z)$ is a constant.
\end{proposition}

{\bf Proof}.
Let $S_1$ be the common domain of definition of
the operators $\widehat X^a\widehat Y^b$, $a \leq 0, b \geq 0$.
\begin{lemma} \label{per2} The space $S_1$ consists of the functions $f(x)$ in $L_2$ which admit an analytic continuation to the upper half plane $y>0$, and decay faster then
$e^{ax}$ for any $a>0$  on each line $x+iy$.
\end{lemma}

{\bf Proof}.
Indeed, it is invariant under multiplication by $e^{bx}$, $b>0$, and
shift by $2\pi ia$, $a>0$, which means that Fourier transform of a function from $S_1$
is invariant under multiplication by $e^{ax}$, $a>0$.
The lemma is proved.
\vskip 3mm

Since $K^5S \subset S$, it follows that $F(x)w \in S\subset S_1$ for
$w \in W$,
and hence
$F(z)$ is analytic in the half plane
$y >0$.
The operator of multiplication by $F(z)$
commutes with the shifts by $2 \pi  i$ and $2 \pi i \hbar$. Thus
it commutes with the shift by $2 \pi  i (n + m\hbar)$, $m,n>0$.
This implies that $F(z)$ is invariant under the shifts by $2\pi i(m+n\hbar)$
where $m+n\hbar >0$.
Thus $F(z)$ is
a constant when $\hbar$ is irrational.
Since $K^5$ depends continuously on $\hbar$, we get
Proposition \ref{per1}, and hence  Theorem \ref{qp4*}.
\vskip 3mm

\section{Algebraic geometry of $  {\cal M}_{0,5}^{\rm cyc}$}

Recall the cross-ratio $r(x_1, x_2, x_3, x_4)$ of four points on ${\Bbb P}^1$:
$$
r(x_1, x_2, x_3, x_4):=
\frac{(x_1 - x_2)(x_3 - x_4)}{(x_1 - x_4)(x_2 - x_3)}, \qquad r(\infty, -1, 0, z) = z.
$$
 It satisfies the relations
$
r(x_1, x_2, x_3, x_4)= r(x_2, x_3, x_4, x_1)^{-1} = -1- r(x_1, x_3, x_2, x_4).
$

Let $ {\cal M}_{0,5}$ be the moduli space of
configurations of five distinct points  on ${\Bbb P}^1$ considered modulo the
action of $PGL_2$. The moduli space $\overline {\cal M}_{0,5}$ is a smooth
algebraic  surface compactifying  $ {\cal M}_{0,5}$.
There are $10$ projective lines $D_{ij}$,  $1 \leq i < j \leq 5$, inside of
$\overline {\cal M}_{0,5}$, forming ``the divisor at infinity'' $D = \cup D_{ij}$.
The line $D_{ij}$ parametrizes configurations of points $(x_1, x_2, x_3, x_4, x_5)$
where ``$x_i$ collides with $x_j$''. So  $\overline {\cal M}_{0,5}-D =  {\cal M}_{0,5}$.

We picture points $x_1, \ldots, x_5$  at the vertexes of an oriented
 pentagon, whose orientation agrees
with the cyclic order of the points. Given a triangulation $T$ of the pentagon,
let us define a pair of rational functions on the surface $\overline {\cal M}_{0,5}$,
assigned to the diagonals of the triangulation.
Given a diagonal $E$, let $z_1, z_2, z_3, z_4$ be the configuration of four points at the vertexes of the rectangle containing $E$ as a diagonal, so that $z_1$ is a vertex of $E$. Then
we set
$$
X_E^T:=  r(z_1, z_2, z_3, z_4).
$$

{\bf Example}. Given a configuration $(\infty, -1, 0, x, y)$, and taking the triangulation
related to the vertex at $\infty$, we get functions
$X = x, Y = (y-x)/x$, see Fig. \ref{qp4*}.

\begin{figure}[ht]
\centerline{\epsfbox{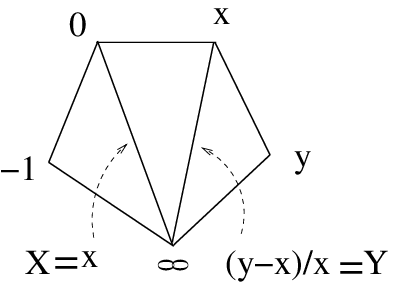}}
\caption{}
\label{qp4}
\end{figure}

\begin{definition} \label{pq1}
$
  {\cal M}^{\rm cyc}_{0,5}:=
\overline {\cal M}_{0,5}-\cup_{c=1}^5D_{c, c+1}$, where $c$ is  modulo $5$.
\end{definition}

The space $  {\cal M}^{\rm cyc}_{0,5}$ is determined by a choice of cyclic order of
configurations of points $(x_1, \ldots, x_5)$.

Let us define embeddings
$
\psi_c: \C^* \times \C^* \hra \overline {\cal M}_{0,5}
$  for $c \in \{1, \ldots, 5\}$.  Set
$$
\psi_1: (X, Y)\lms (\infty, -1, 0, X, X(1+Y)).
$$
One easily checks that it is an embedding.
The map $\psi_c$ is obtained from $\psi_1$ by the cyclic shift of the configuration of five points
by $2c$. So it is also an embedding.

\begin{figure}[ht]
\centerline{\epsfbox{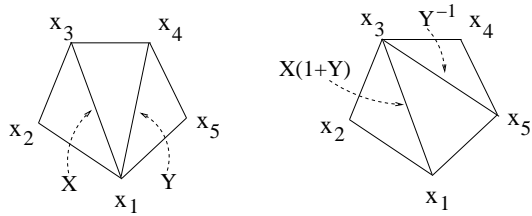}}
\caption{Change of the coordinates under a flip.}
\label{qp3}
\end{figure}

\vskip 3mm
The following function is regular on the surface $  {\cal M}^{\rm cyc}_{0,5}$:
$$
X_{a,b; c}:= r(x_{c}, x_{c+1}, x_{c+2}, x_{c+3} )^ar(x_{c}, x_{c+2}, x_{c+3}, x_{c+4} )^b,
\qquad a\geq 0, b\leq 0.
$$
Indeed, the poles of the first factor are at the divisor $D_{c, c+1}\cup D_{c+2, c+3}$, and
the poles of the second one are
 at the divisor $D_{c+2, c+3}\cup D_{c, c+4}$. The set of functions $\{X_{a,b; c}\}$
coincides with the one defined in Section 2.2. Indeed, one checks
 this  for $c=1$ using Fig. \ref{qp4},
and use equivariance with respect to the shifts and Fig. \ref{qp3}.

\begin{theorem} \label{pq2} (i) The surface $  {\cal M}^{\rm cyc}_{0,5}$
is the union of the five open subsets $\psi_c(\C^* \times \C^*)$
in  $\overline {\cal M}_{0,5}$.

(ii) The functions $X_{a,b; c}$, where $a, b\in Z, a\geq 0, b \leq 0$ and $c$ is mod $5$ form a basis of the space of regular functions on the surface $  {\cal M}^{\rm cyc}_{0,5}$.
\end{theorem}

{\bf Proof}. (i) Straightforward.

\begin{figure}[ht]
\centerline{\epsfbox{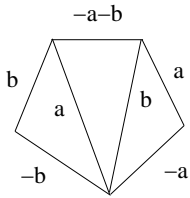}}
\caption{The weight diagram of a basis monomial; $a,b \geq 0$.}
\label{qp5*}
\end{figure}

(ii) The algebra of regular functions on ${\cal M}_{0,5}$ is defined as follows.
Take the configuration space ${\rm Conf}_5(V_2)$ of $5$-tuples of vectors
$(v_1, \ldots, v_5)$ in generic position in a two-dimensional symplectic vector space $V_2$, modulo
the $SL_2$-action. The group $H:= (\C^*)^5$ acts on it by multiplying each
vector $v_i$ by a number $\lambda_i$.
Let $\Delta_{ij}$ is the area in $V_2$ of the parallelogram $\langle v_i, v_j\rangle$. Then
\begin{equation} \label{pop}
\Z[{\cal M}_{0,5}] = \Z[{\rm Conf}_5(V_2)]^H = \Z[\Delta^{\pm 1}_{ij}]^H.
\end{equation}
 The subspace $\Z[  {\cal M}^{\rm cyc}_{0,5}]$  is spanned by the monomials
\begin{equation} \label{pop1}
\prod_{1\leq i<j\leq 5}\Delta_{ij}^{a_{ij}},
\end{equation}
where
$a_{i, i+1} \in \Z$, and $ a_{ij} \in \Z_{\geq 0}$ unless $j = i\pm 1$ mod $5$.
Write the integers $a_{ij}$ on the diagonals and sides of the pentagon.
Call them the {\it weights} and the corresponding picture the {\it chord diagram}.
The $H$-invariance means
 that the sum of the weights assigned to the edges and sides sharing a vertex is $0$.
We erase diagonals of  weight $0$.
A monomial is {\it regular} if its chord diagram has no intersecting diagonals.
Using the Pl\"ucker relations $\Delta_{ac}\Delta_{bd} = \Delta_{ab}\Delta_{cd} +
\Delta_{ad}\Delta_{bc}$, $1 \leq a<b<c<d\leq 5$,
and arguing by induction on the sum of the products of the weights of
diagonals in the intersection points, we reduce any sum of monomials
(\ref{pop1}) to a sum of the regular ones. An easy argument with the
``sum of the weights at a vertex equals zero'' equations shows that
for a regular monomial there exists a vertex of the pentagon such that its
 weights are as on Fig. \ref{qp5*}. So the functions $X_{a,b; c}$
span the space of regular functions on $  {\cal M}^{\rm cyc}_{0,5}$. To check that they
are linearly independent, look at the monomials with maximum value of
$a+b$. The theorem is proved.

\vskip 3mm
{\bf The quantized ${\cal M}^{\rm cyc}_{0,5}$ at roots of unity.}
Assume that $q$ is a primitive $N$-th root of unity. Then the functions
$x_{a,b; c}:= X_{Na,Nb; c}$ generate the center of the algebra ${\Bbb L}_q$.
In particular
$x:= X^N, y:=Y^N$ are in the center. One checks (\cite{FG2}, Section 3) that
the elements $x,y$ behave under flips just like the corresponding coordinates on
${\cal M}^{\rm cyc}_{0,5}$.
Therefore the spectrum of the center 
of ${\Bbb L}_q$ is identified with ${\cal M}^{\rm cyc}_{0,5}$.
Restricting to an affine chart of ${\cal M}^{\rm cyc}_{0,5}$ with coordinates $(\alpha, \beta)$
we see that the localization of the algebra ${\Bbb L}_q$ at this chart is identified
with the  algebra generated by $X,Y$ with the relations
$X^N=\alpha, Y^N = \beta, XY = q^2YX$. It is well know that it is a sheaf of central simple
algebras over $\C^* \times \C^*$. So we get 

\begin{proposition}
Let $q$ be a root of unity. Then the algebra ${\Bbb L}_q$ gives rise to 
 a sheaf of Azumaya algebras on ${\cal M}^{\rm cyc}_{0,5}$.
\end{proposition}

The real positive part of ${\cal M}^{\rm cyc}_{0,5}$ is given by 
configurations of points $(\infty, -1, 0, x, y)$ with $0 < x < y$. Its 
closure 
in $\overline {\cal M}_{0,5}(\R)$ is the pentagon. Its sides are real segments on the divisors 
$D_{c, c+1}$. 

%One can recover the  moduli space $\overline {\cal M}_{0,5}$ from the algebra 
%${\Bbb L}$ as follows. For each $c = 1, \ldots , 5$ consider the subalgebra 
%${\Bbb L}^{(c)}$ in ${\Bbb L}$ given by the functions regular 
%at all the the divisors $D_{c', c'+1}$ different from $D_{c, c+1}$. 

\vskip 3mm
Address: A.B. Goncharov, 
Dept of Mathematics, Brown University, Providence RI 02906, USA. e-mail sasha@math.brown.edu.

\end{document}